\documentclass{article}
\usepackage{latexsym}
\usepackage{amsmath}
\usepackage{amssymb}
\usepackage{epsfig}
\usepackage{amsfonts}

\begin{document}
\title{A conjecture in the problem of rational definite summation}
\author{Mark van Hoeij}
\maketitle
\begin{abstract}
A conjecture is given that, if true, could lead to an
algorithm for computing definite sums of rational functions.
\end{abstract}
\section{The conjecture}

Let $C$ be a field of characteristic 0,
not necessarily algebraically closed.
Let $n, k$ be variables, and consider the ring $C(n)[k]$
of polynomials in $k$, whose coefficients are rational
functions in $n$.

For all $a,b,c \in \mathbf{Q}$ with $a \neq 0$ we define
the $C(n)$-algebra homomorphism
\[ \psi_{a,b,c}: C(n)[k] \rightarrow C(n)[k] \]
as follows:
\[ \psi_{a,b,c}(k) = ak + bn + c. \]

If $B_1, B_2$ are irreducible elements of $C(n)[k]$,
then we define the following equivalence relation:
$B_1 \sim_p B_2$ if there exist
$a,b,c \in \mathbf{Q}$, with $a \neq 0$ such that
$\psi_{a,b,c}(B_1) = s B_2$ for some unit $s$ (i.e.
$s$ is a nonzero element of $C(n)$).

Let $F \in C(n,k)$ be a nonzero rational function in $n$
and $k$. The {\em partial fraction decomposition} of
$F = F(n,k)$ over $C(n)$ has the following form:
\[ F = Q + F_1 + \cdots + F_t \]
where $Q \in C(n)[k]$, and each $F_i$ is a {\em
partial fraction term} which means the following:
there exist a positive integer $d_i$, non-zero
$A_i,B_i \in C(n)[k]$ with ${\rm degree}_k(A_i) < {\rm
degree}_k(B_i)$,
with $B_i$ irreducible, such that
\[ F_i = A_i/{B_i}^{d_i}. \]
We call two partial-fraction-terms $F_i$ and $F_j$
{\em equivalent} $F_i \sim F_j$ when
\[ d_i = d_j {\rm \ \ and \ \ } B_i \sim_p B_j. \]
Suppose that there are $u$ distinct equivalence
classes among $F_1,\ldots,F_t$, and number these
equivalence classes $1,\ldots,u$. Let $G_i$ be the
sum of all $F_j$ in the $i$'th equivalence class.
So we can write
\[ F = Q + G_1 + \cdots G_u. \]

Now $\sum_{k=0}^n F(n,k)$ is a function in just one
variable $n$. This function is
only defined for nonnegative integers $n$ for which
the denominator of $F(n,k)$ does not vanish for any
$k \in \{0,1,\ldots,n\}$.
Now assume that this is the case for $n \gg 0$.

The question now is: Under this assumption,
how to decide if this function $\sum_{k=0}^n F$
is a rational function?
More precisely: is there a rational function $R(n) \in C(n)$
that takes the same values as $\sum_{k=0}^n F$
whenever both are defined? \\[10pt]
{\bf Conjecture:} If $\sum_{k=0}^n F$ is a rational function,
then so is $\sum_{k=0}^n G_i$ for every
$i \in \{1,\ldots,u\}$. \\[10pt]

The conjecture says that the problem of deciding if
$\sum_{k=0}^n F$ is a rational function (and if so,
finding that rational function) can be reduced
to the same question for the $\sum_{k=0}^n G_i$.
This reduces the rational summation problem to a set of smaller
summation problems  $\sum_{k=0}^n G_i$, each of which
involves only one equivalence class of partial-fraction-terms.

\section{A few remarks}
Let $B$ be an irreducible element of $C(n)[k]$.
Define ${\rm Aut}(B)$ as the set of all $\psi_{a,b,c}$
with $a,b,c \in \mathbf{Q}$, $a \neq 0$, such that
$\psi_{a,b,c}(B) = s B$ for some unit $s$.
It is easy to see that ${\rm Aut}(B)$ is a group under
composition, and that the map
\[ \pi:  {\rm Aut}(B) \rightarrow \mathbf{Q}^* \]
sending $\psi_{a,b,c}$ to $a$ is an injective group
homomorphism. From this it follows that ${\rm Aut}(B)$
has either 1, 2, or $\infty$ many elements.
We call $B$ of {\em generic type} if ${\rm Aut}(B)$ has 1 element,
{\em symmetric type} if ${\rm Aut}(B)$ has 2 elements,
and {\em rational type} if ${\rm Aut}(B)$ has $\infty$
many elements.
Note that $B$ is of rational type if and only if $B \sim_p k$.

$G_i$ is a sum of equivalent terms $A_j / {B_j}^{d_j}$.
We call $G_i$ of generic, symmetric, resp. rational type
if the $B_j$ appearing in $G_i$ are of generic, symmetric, resp. rational
type.
In order to calculate $\sum_{k=0}^n G_i$ it helps
to distinguish these three types.  The rational type allows
more kinds of cancelations among equivalent terms $A_j / {B_j}^{d_j}$
than the symmetric type, which
in turn allows more kinds of cancelations than the generic type.
We give an example of cancelation for the generic type:
Let $B(n,k)$ be an arbitrary irreducible
polynomial. For example, $B(n,k) = k^3 + kn + 1$. Consider:
\begin{equation}
 \label{1}
 \frac{1}{B(n, k)} - \frac{1}{B(n, n-k)}
\end{equation}
or
\begin{equation}
 \label{2}
 \frac{1}{B(n, k)} + \frac{1}{B(n,k+n+1)} -
 \frac{1}{B(n, 2k)} - \frac{1}{B(n, 2k + 1)}.
\end{equation}
These are non-zero rational functions, whose $\sum_{k=0}^n$
is $0$. The algorithm in~\cite{H} would detect
cancelations like the one in equation~(\ref{1}) but not
cancelations like the one in~(\ref{2}).  This problem was
the starting point for this paper.

The rational and symmetric types allow additional kinds of cancelation.
Here is an example for the rational type:
\[
	1/(k+1) + 2/(k+2+n) - 2/(2k+1)
\]
has sum 0 because it can be rewritten as
\[
	2/(k+1) + 2/(k+2+n) - 2/(2k+1) - 2/(2k+2).
\]
We can give examples where the sum is a non-zero rational
function by modifying
the above examples in such a way that all but a constant
number of terms cancel. For example, the sum of:
\begin{equation}
	1/(k+3) + 2/(k+2+n) - 2/(2k+1)
\label{nonz}
\end{equation}
differs from the sum of $1/(k+1) + 2/(k+2+n) - 2/(2k+1)$
(which was 0) by $D = 1/(n+3) + 1/(n-1+3) - 1/(0+1) - 1/(1+1)$.
To develop an algorithm for rational definite summation, one must
identify all possible cancelations and write procedures
for each of them. Then, for inputs such
as~(\ref{nonz}), the algorithm must
determine if this input can be modified into something that has sum 0,
and if so, calculate the difference $D$ like in the example.
One could then try to prove that the algorithm
completely solves the problem, that it determines for $F \in C(n,k)$
if $\sum_{k=0}^n F$ is a rational function or not, and if so, finds
that rational function.
For a preliminary implementation of this algorithm, see: \\
\ \ \ http://www.math.fsu.edu/$\tilde{\mbox{ \ }}$hoeij/SumRat/ \\
It implements several kinds of cancelations. It is not clear
(even if one assumes the conjecture)
if this algorithm is complete or not. \\

\noindent {\bf A second question:}
Suppose that the denominator of $F$ has no irreducible factors $B \in
C[n,k]$ for which there exists $P \in C[k]$ with $B \sim_p P$.
Suppose furthermore that $R(n) := \sum_{k=0}^n F(n,k)$
satisfies a linear recurrence relation $\sum_{i=0}^m a_i(n) R(n+i) = 0$
for some polynomials $a_i(n) \in C[n]$, not all zero. Does it then follow
that $R(n)$ is a rational function?

\section{An example}
Let
\[ B := k + n(n+1)/2, \ \ \ \ F := 1/B, \ \ \ \
R(n) := \sum_{k=0}^n F \ \ \ \ {\rm where \ }n>0.\]
Then $R(n)$ does not satisfy a linear recurrence with polynomial 
coefficients. \\

\noindent {\bf Proof:}
Suppose that $m$ is a positive integer
and that $a_0(n), \ldots, a_m(n)$ are in ${\mathbf Z}[n]$
such that $\sum_{i=0}^m a_i(n) R(n+i)$ is identically 0.

Let $S_n$ be the set $\{n(n+1)/2, 1 + n(n+1)/2, \ldots, n + 
n(n+1)/2\}$,
so $R(n) = \sum_{i \in S_n} 1/i$.
Let $A_n,B_n$ be positive integers, with $A_n/B_n = R(n)$
and ${\rm gcd}(A_n,B_n)=1$.
The sets $S_1,S_2,\ldots$ are disjoint and their union is the set of
all positive integers.
Define $p_n$ as the product of all prime numbers in $S_n$.
For $0<i<n$,
these primes are larger than ${\rm max}(S_i)$, so
${\rm gcd}(p_n, B_i)=1$.

Assume that $a_0(n),\ldots,a_m(n) \in {\mathbf Z}[n]$,
and that $a_m(n)$ is not identically zero.
Then there exists an integer $n_1$ such
that $a_m(n) \neq 0$ for all $n>n_1$.
If $\sum_{i=0}^m a_i(n) R(n+i) = 0$  then
$a_m(n)$ must be divisible by $p_{m+n}$ because the factors of
$p_{m+n}$ appear in the denominator of $R(m+n)$ but not in
the denominators of $R(n), R(n+1), \ldots, R(m-1+n)$.
Hence $p_{m+n} \leq | a_m(n) |$.
Since each prime factor in $p_{m+n}$ is $> \frac12 n^2$
we see that the number of primes in $S_{m+n}$ is
at most ${\rm ln}( | a_m(n) | ) / {\rm ln}(\frac12 n^2)$,
which tends to ${\rm degree}(a_m)/2$ for large $n$.
This means that the density of primes in
$S_n$ is much less than $1/{\rm ln}({\rm max}(S_n))$,
contradicting the prime number theorem (note that
$S_1,S_2,\ldots$ are disjoint and their union is the set of
all positive integers).

One sees by induction that $a_{m-1},\ldots,a_0$ must
be identically 0 as well. So $R(n)$ does not satisfy any
nonzero linear recurrence relation over ${\mathbf Z}[n]$.
There can not exist a recurrence relation over ${\mathbf C}[n]$
either, because if $K \subset {\mathbf C}$ is the coefficient
field of such a recurrence relation, then one can reduce the
transcendence degree of the finitely generated field $K$ by switching
to a residue field of a valuation.
After finitely many steps, we may assume that the transcendence degree 
is 0, so $K$ is a number field,
Then one can multiply the recurrence by an element of $K$ in order
to create at least one polynomial $a_i(n)$ with leading coefficient 1.
Then, after taking the trace over the field of rational numbers,
one obtains a nonzero recurrence over ${\mathbf Q}[n]$, and by
multiplying out the denominators one obtains a recurrence over
${\mathbf Z}[n]$. \\
Hopefully primes (or prime ideals) could be used for the conjecture
as well. \\

\noindent{\bf Acknowledgments:}
I would like to thank Sergei Abramov, Ha Le, Sergei Tsarev and Eugene Zima
for discussions on these topics.

\end{document}